\documentclass[10pt, a4paper]{article}

\usepackage{amsfonts, latexsym, amsmath, amssymb, fancyhdr, amsthm, graphicx, subfigure}
\usepackage[english]{babel}

\newtheorem{theorem}{Theorem}[section]
\newtheorem{lemma}{Lemma}[section]
\newtheorem{proposition}{Proposition}[section]

\newtheorem{remark}{Remark}

\title{\textbf{Global Asymptotic Stability and Hopf Bifurcation for a Blood Cell Production Model}\thanks{Published in Mathematical Biosciences and Engineering, vol 3, iss 2, 325-346 (2006)}}
\author{Fabien {\sc Crauste}}
\date{Year 2005}

\begin{document}

\maketitle

\begin{center}
\emph{Laboratoire de Math\'ematiques Appliqu\'ees, UMR 5142,} \\
\emph{Universit\'e de Pau et des Pays de l'Adour,}\\
\emph{Avenue de l'universit\'e, 64000 Pau, France.}\\
\emph{ANUBIS project, INRIA--Futurs}\\
\emph{E-mail: fabien.crauste@univ-pau.fr}\\
\end{center}

\quad

\begin{abstract}
We analyze the asymptotic stability of a nonlinear system of two differential equations with delay
describing the dynamics of blood cell production. This process takes place in the bone marrow where
stem cells differentiate throughout divisions in blood cells. Taking into account an explicit role
of the total population of hematopoietic stem cells on the introduction of cells in cycle, we are
lead to study a characteristic equation with delay-dependent coefficients. We determine a necessary
and sufficient condition for the global stability of the first steady state of our model, which
describes the population's dying out, and we obtain the existence of a Hopf bifurcation for the
only nontrivial positive steady state, leading to the existence of periodic solutions. These latter
are related to dynamical diseases affecting blood cells known for their cyclic nature.
\end{abstract}

\bigskip{}


\section{Introduction}

Blood cell production process is based upon the differentiation of so-called hematopoietic stem
cells, located in the bone marrow. These undifferentiated and unobservable cells have unique
capacities of differentiation (the ability to produce cells committed to one of the three blood
cell types: red blood cells, white cells or platelets) and self-renewal (the ability to produce
cells with the same properties).

Mathematical modelling of hematopoietic stem cells dynamics has been introduced at the end of the
seventies by Mackey \cite{m1978}. He proposed a system of two differential equations with delay
where the time delay describes the cell cycle duration. In this model, hematopoietic stem cells are
separated in proliferating and nonproliferating cells, these latter being introduced in the
proliferating phase with a nonlinear rate depending only upon the nonproliferating cell population.
The resulting system of delay differential equations is then uncoupled, with the nonproliferating
cells equation containing the whole information about the dynamics of the hematopoietic stem cell
population. The stability analysis of the model in \cite{m1978} highlighted the existence of
periodic solutions, through a Hopf bifurcation, describing in some cases diseases affecting blood
cells, characterized by periodic oscillations \cite{hdm1998}.

The model of Mackey \cite{m1978} has been studied by many authors, mainly since the beginning of
the nineties. Mackey and Rey \cite{mr1993, mr1995_2, mr1995_1} numerically studied the behavior of
a structured model based on the model in \cite{m1978}, stressing the existence of strange behaviors
of the cell populations (like oscillations, or chaos). Mackey and Rudnicky \cite{mr1994, mr1999}
developed the description of blood cell dynamics through an age-maturity structured model,
stressing the influence of hematopoietic stem cells on blood production. Their model has been
further developed by Dyson et al. \cite{dvw2000_1, dvw2000_2, dvw2002}, Adimy and Pujo-Menjouet
\cite{ap2003}, Adimy and Crauste \cite{ac2003_2, ac2005} and Adimy et al. \cite{acp2005}. Recently,
Adimy et al. \cite{acr2005_2, acr2005} studied the model proposed in \cite{m1978} taking into
account that cells in cycle divide according to a density function (usually gamma distributions
play an important role in cell cycles durations), contrary to what has been assumed in the
above-cited works, where the division has always been assumed to occur at the same time.

More recently, Pujo-Menjouet and Mackey \cite{pm2004} and Pujo-Menjouet et al. \cite{pbm2005} gave
a better insight into the model of Mackey \cite{m1978}, highlighting the role of each parameter of
the model on the appearance of oscillations and, more particularly, of periodic solutions, when the
model is applied to the study of chronic myelogenous leukemia \cite{fm1999}.

Contrary to the assumption used in all of the above-cited works, we study, in this paper, the model
introduced by Mackey \cite{m1978} considering that the rate of introduction in the proliferating
phase, which contains the nonlinearity of this model, depends upon the total population of
hematopoietic stem cells, and not only upon the nonproliferating cell population. The introduction
in cell cycle is partly known to be a consequence of an activation of hematopoietic stem cells due
to molecules fixing on them. Hence, the entire population is in contact with these molecules and it
is reasonable to think that the total number of hematopoietic stem cells plays a role in the
introduction of nonproliferating cells in the proliferating phase.

The first consequence is that the model is not uncoupled, and the nonproliferating cell population
equation does not contain the whole information about the dynamics of blood cell production,
contrary to the model in \cite{m1978, pbm2005, pm2004}. Therefore, we are lead to the study of a
modified system of delay differential equations (system (\ref{eqS})--(\ref{eqN2})), where the delay
describes the cell cycle duration, with a nonlinear part depending on one of the two populations.

Secondly, while studying the local asymptotic stability of the steady states of our model, we have
to determine roots of a characteristic equation taking the form of a first degree exponential
polynomial with delay-dependent coefficients. For such equations, Beretta and Kuang \cite{bk2002}
developed a very useful and powerful technic, that we will apply to our model.

Our aim is to show, through the study of the steady states' stability, that our model, described in
(\ref{eqS})--(\ref{eqN2}), exhibits similar properties than the model in \cite{m1978} and that it
can be used to model blood cells production dynamics with good results, in particularly when one is
interested in the appearance of periodic solutions in blood cell dynamics models. We want to point
out that the usually accepted assumption about the introduction rate may be limitative and that our
model can display interesting dynamics, such as stability switches, that have never been noted
before.

The present work is organized as follows. In the next section we present our model, stated in
equations (\ref{eqS}) and (\ref{eqN2}). We then determine the steady states of this model. In
section \ref{sceas}, we linearize the system (\ref{eqS})--(\ref{eqN2}) about a steady state and we
deduce the associated characteristic equation. In section \ref{stss}, we establish necessary and
sufficient conditions for the global asymptotic stability of the trivial steady state (which
describes the extinction of the hematopoietic stem cell population). In section \ref{spss}, we
focus on the asymptotic stability of the unique nontrivial steady state. By studying the existence
of pure imaginary roots of a first degree exponential polynomial with delay-dependent coefficients,
we obtain the existence of a critical value of the time delay for which a Hopf bifurcation occurs
at the positive steady state, leading to the appearance of periodic solutions. Using numerical
illustrations, we show how these solutions can be related to periodic hematological diseases in
section \ref{snum}, and we note the existence of a stability switch. We conclude with a discussion.

\section{A Nonlinear Model of Blood Cell Production}

Let consider a population of hematopoietic stem cells, located in the bone marrow. These cells
actually perform a succession of cell cycles, in order to differentiate in blood cells (white
cells, red blood cells and platelets). According to early works, by Burns and Tannock \cite{bt} for
example, we assume that cells in cycle are divided in two groups: proliferating and
nonproliferating cells. The respective proliferating and nonproliferating cell populations are
denoted by $P$ and $N$.


All hematopoietic stem cells die with constant rates, namely $\gamma>0$ for proliferating cells and
$\delta>0$ for nonproliferating cells. These latter are introduced in the proliferating phase, in
order to mature and divide, with a rate $\beta$. At the end of the proliferating phase, cells
divide in two daughter cells which immediately enter the nonproliferating phase.

Then the populations $P$ and $N$ satisfy the following evolution equations (see Mackey \cite{m1978}
or Pujo-Menjouet and Mackey \cite{pm2004}),
\begin{eqnarray}
\displaystyle\frac{dP}{dt}(t)&=&-\gamma P(t) +\beta N(t) -e^{-\gamma\tau}\beta
N(t-\tau),\label{eqP}\\
\displaystyle\frac{dN}{dt}(t)&=&-\delta N(t) -\beta N(t) +2e^{-\gamma\tau}\beta
N(t-\tau).\label{eqN}
\end{eqnarray}
In each of the above equations, $\tau$ denotes the average duration of the proliferating phase. The
term $e^{-\gamma\tau}$ then describes the survival rate of proliferating cells. The last terms in
the right hand side of equations (\ref{eqP}) and (\ref{eqN}) account for cells that have performed
a whole cell cycle and leave (enter, respectively) the proliferating phase (the nonproliferating
phase, respectively). These cells are in fact nonproliferating cells introduced in the
proliferating phase a time $\tau$ earlier. The factor 2 in equation (\ref{eqN}) represents the
division of each proliferating hematopoietic stem cell in two daughter cells.

We assume that the rate of introduction $\beta$ depends upon the total population of hematopoietic
stem cells, that we denote by $S$. With our notations, $S=P+N$. This assumption stresses the fact
that the nature of the trigger signal for introduction in the proliferating phase is the result of
an action on the entire cell population. For example, it can be caused by molecules entering the
bone marrow and fixing on hematopoietic stem cells, activating or inhibiting their proliferating
capacity. This occurs in particularly for the production of red blood cells. Their regulation is
mainly mediated by an hormone (a growth factor, in fact) called erythropoietin, produced by the
kidneys under a stimulation by circulating blood cells (see B\'elair et al. \cite{bmm1995}, Mahaffy
et al. \cite{mbm1998}).

Hence we assume that
\begin{displaymath}
\beta=\beta(S(t)).
\end{displaymath}
The function $\beta$ is supposed to be continuous and positive on $[0,+\infty)$, and strictly
decreasing. This latter assumption describes the fact that the less hematopoietic stem cells in the
bone marrow, the more cells introduced in the proliferative compartment \cite{m1978, pm2004}.
Furthermore, we assume that
\begin{displaymath}
\lim_{S\to\infty}\beta(S)=0.
\end{displaymath}

Adding equations (\ref{eqP}) and (\ref{eqN}) we can then deduce an equation satisfied by the total
population of hematopoietic stem cells $S(t)$. We assume, for the sake of simplicity, that
proliferating and nonproliferating cells die with the same rate, that is $\delta=\gamma$. Then the
populations $N$ and $S$ satisfy the following nonlinear system with time delay $\tau$,
corresponding to the cell cycle duration,
\begin{eqnarray}
\displaystyle\frac{dS}{dt}(t)&=&-\delta S(t) +e^{-\delta\tau}\beta(S(t-\tau))N(t-\tau), \label{eqS}
\vspace{1ex}\\
\displaystyle\frac{dN}{dt}(t)&=&-\delta N(t) -\beta(S(t))N(t)
+2e^{-\delta\tau}\beta(S(t-\tau))N(t-\tau).\label{eqN2}
\end{eqnarray}

From Hale and Verduyn lunel \cite{haleverduynlunel}, for each continuous initial condition, system
(\ref{eqS})--(\ref{eqN2}) has a unique continuous solution $(S(t), N(t))$, well-defined for
$t\geq0$.

\begin{lemma}\label{lempos}
For all nonnegative initial condition, the unique solution $(S(t), N(t))$ of
(\ref{eqS})--(\ref{eqN2}) is nonnegative.
\end{lemma}

\begin{proof}
First assume that there exists $\xi>0$ such that $N(\xi)=0$ and $N(t)>0$ for $t<\xi$. Then, from
(\ref{eqN2}) and since $\beta$ is a positive function,
\begin{displaymath}
\frac{dN}{dt}(\xi)=2e^{-\delta\tau}\beta(S(\xi-\tau))N(\xi-\tau)>0.
\end{displaymath}
Consequently, $N(t)\geq0$ for $t>0$.

If there exists $\zeta>0$ such that $S(\zeta)=0$ and $S(t)>0$ for $t<\zeta$, then the same
reasoning, using (\ref{eqS}), leads to
\begin{displaymath}
\frac{dS}{dt}(\zeta)=e^{-\delta\tau}\beta(S(\zeta-\tau))N(\zeta-\tau)>0,
\end{displaymath}
and we deduce that $S(t)\geq0$ for $t>0$.
\end{proof}

\begin{remark}
The positivity of $S$ and $N$, solutions of system (\ref{eqS})--(\ref{eqN2}), does not a priori
implies that $P=S-N$ is nonnegative.

Using a classical variation of constant formula, the solutions $P(t)$ of (\ref{eqP}) are given, for
$t\geq0$, by
\begin{displaymath}
P(t)=e^{-\delta t}P(0)+e^{-\delta t}\int_0^t
e^{\delta\theta}\beta(S(\theta))N(\theta)-e^{\delta(\theta-\tau)}\beta(S(\theta-\tau))N(\theta-\tau)
d\theta.
\end{displaymath}
Setting the change of variable $\sigma=\theta-\tau$, we obtain
\begin{displaymath}
P(t)=e^{-\delta t}\left[P(0)-\int_{-\tau}^0
e^{\delta\theta}\beta(S(\theta))N(\theta)d\theta\right]+e^{-\delta t}\int_{t-\tau}^t
e^{\delta\theta}\beta(S(\theta))N(\theta) d\theta.
\end{displaymath}
Consequently, $P(t)\geq0$ for $t\geq0$ if
\begin{displaymath}
P(0)\geq \int_{-\tau}^0 e^{\delta\theta}\beta(S(\theta))N(\theta)d\theta,
\end{displaymath}
that is if
\begin{displaymath}
S(0)\geq N(0)+\int_{-\tau}^0 e^{\delta\theta}\beta(S(\theta))N(\theta)d\theta.
\end{displaymath}
This condition is biologically relevant since $\int_{-\tau}^0
e^{\delta\theta}\beta(S(\theta))N(\theta)d\theta$ represents the population of cells that have been
introduced in the proliferating phase at time $\theta\in[-\tau,0]$ and that have survived at time
$t=0$. Hence, from a biological point of view, the population in the proliferating phase at time
$t=0$ should be larger than this quantity.
\end{remark}

In the stability analysis of system (\ref{eqS})--(\ref{eqN2}), the existence of stationary
solutions, called steady states, is relevant since these particular solutions are potential limits
of system (\ref{eqS})--(\ref{eqN2}).

A steady state of system (\ref{eqS})--(\ref{eqN2}) is a solution $(\overline{S},\overline{N})$
satisfying
\begin{displaymath}
\frac{d\overline{S}}{dt}=\frac{d\overline{N}}{dt}=0.
\end{displaymath}
Let $(\overline{S},\overline{N})$ be a steady state of (\ref{eqS})--(\ref{eqN2}). Then
\begin{eqnarray}
\delta\overline{S}&=&e^{-\delta\tau}\beta(\overline{S})\overline{N},\label{ss1}\\
(\delta+\beta(\overline{S}))\overline{N}&=&2e^{-\delta\tau}\beta(\overline{S})\overline{N}.
\label{ss2}
\end{eqnarray}

We immediately notice that $(0,0)$ is a steady state of system (\ref{eqS})--(\ref{eqN2}), that we
will denote, in the following, by $E^0$. This steady state always exists. It describes the
extinction of the hematopoietic stem cell population.

Assume that $(\overline{S},\overline{N})$ is a steady state of (\ref{eqS})--(\ref{eqN2}) with
$\overline{S}, \overline{N}\neq0$. Then, from (\ref{ss1}) and (\ref{ss2}),
\begin{displaymath}
(2e^{-\delta\tau}-1)\beta(\overline{S})=\delta \qquad \textrm{and} \qquad
\overline{N}=\frac{\delta\overline{S}}{e^{-\delta\tau}\beta(\overline{S})}.
\end{displaymath}
A necessary condition to obtain a nontrivial steady state is then $2e^{-\delta\tau}-1>0$, that is
\begin{displaymath}
\tau<\frac{\ln(2)}{\delta}.
\end{displaymath}
Under this condition, since the function $\beta$ is decreasing, positive, and tends to zero at
infinity, there exists $\overline{S}>0$ satisfying
\begin{equation}\label{eqss}
(2e^{-\delta\tau}-1)\beta(\overline{S})=\delta,
\end{equation}
if and only if
\begin{equation}\label{condexist}
(2e^{-\delta\tau}-1)\beta(0)>\delta.
\end{equation}
One can easily check that condition (\ref{condexist}) is equivalent to
\begin{equation}\label{condexist2}
\delta<\beta(0) \qquad \textrm{and} \qquad
0\leq\tau<\overline{\tau}:=\frac{1}{\delta}\ln\left(\frac{2\beta(0)}{\delta+\beta(0)}\right).
\end{equation}
In this case, $\overline{S}>0$ solution of (\ref{eqss}) is unique, and
\begin{displaymath}
\overline{N}=\frac{2e^{-\delta\tau}-1}{e^{-\delta\tau}}\overline{S}.
\end{displaymath}
One can check that $\overline{N}\leq\overline{S}$. These results are summed up in the next
proposition.

\begin{proposition}\label{prop1}
If inequality (\ref{condexist}) holds true, the system (\ref{eqS})--(\ref{eqN2}) has exactly two
steady states: $E^0=(0,0)$ and $E^*=(S^*,N^*)$, where $S^*>0$ is the unique solution of equation
(\ref{eqss}) and $N^*=(2e^{-\delta\tau}-1)e^{\delta\tau}S^*$.

If
\begin{displaymath}
(2e^{-\delta\tau}-1)\beta(0)\leq\delta,
\end{displaymath}
then system (\ref{eqS})--(\ref{eqN2}) has only one steady state, namely $E^0=(0,0)$.
\end{proposition}

In the next section, we linearize the system (\ref{eqS})--(\ref{eqN2}) about one of its steady
states in order to analyze its local asymptotic stability.

\section{Linearization and Characteristic Equation}\label{sceas}

We are interested in the asymptotic stability of the steady states of system
(\ref{eqS})--(\ref{eqN2}). To that aim, we linearize system (\ref{eqS})--(\ref{eqN2}) about one of
its steady state and we determine the associated characteristic equation. We assume that $\beta$ is
continuously differentiable on $[0,+\infty)$.

Let $(\overline{S}, \overline{N})$ be a steady state of system (\ref{eqS})--(\ref{eqN2}). From
Proposition \ref{prop1}, $(\overline{S}, \overline{N})$ is either $E^0$ or $E^*$.

The linearization of system (\ref{eqS})--(\ref{eqN2}) about $(\overline{S}, \overline{N})$ leads to
the following system,
\begin{eqnarray}
\displaystyle\frac{dS}{dt}(t)&=&-\delta S(t) +e^{-\delta\tau}\beta(\overline{S})N(t-\tau)+e^{-\delta\tau}\overline{N}\beta^{\prime}(\overline{S})S(t-\tau),\label{eql1}\\
\displaystyle\frac{dN}{dt}(t)&=&-(\delta+\beta(\overline{S}))N(t)
-\overline{N}\beta^{\prime}(\overline{S})S(t)\nonumber\\
&&+2e^{-\delta\tau}\left[\beta(\overline{S})N(t-\tau)+\overline{N}\beta^{\prime}(\overline{S})S(t-\tau)\right],\label{eql2}
\end{eqnarray}
where we have used the notations $S(t)$ and $N(t)$ instead of $S(t)-\overline{S}$ and
$N(t)-\overline{N}$ for the sake of simplicity.

The system (\ref{eql1})--(\ref{eql2}) can be written
\begin{displaymath}
\left(\begin{array}{c}
\displaystyle\frac{dS}{dt}(t)\vspace{1ex}\\
\displaystyle\frac{dN}{dt}(t)
\end{array}\right)
= \mathcal{A}_1 \left(\begin{array}{c}
S(t)\vspace{1ex}\\
N(t)
\end{array}\right)
+ \mathcal{A}_2 \left(\begin{array}{c}
S(t-\tau)\vspace{1ex}\\
N(t-\tau)
\end{array}\right),
\end{displaymath}
where
\begin{displaymath}
\mathcal{A}_1 := \left(\begin{array}{cc}
-\delta & 0 \vspace{1ex}\\
-\alpha & -(\delta+\beta(\overline{S}))
\end{array}\right), \qquad  \mathcal{A}_2 := e^{-\delta\tau}\left(\begin{array}{cc}
\alpha & \beta(\overline{S}) \vspace{1ex}\\
2\alpha & 2\beta(\overline{S})
\end{array}\right),
\end{displaymath}
and
\begin{equation}\label{alpha}
\alpha = \alpha(\overline{N},\overline{S}):=\overline{N}\beta^{\prime}(\overline{S}).
\end{equation}
The characteristic equation of system (\ref{eql1})--(\ref{eql2}) associated with the steady state
$(\overline{S},\overline{N})$ is defined by
\begin{displaymath}
\det(\lambda -\mathcal{A}_1-e^{-\lambda\tau}\mathcal{A}_2)=0.
\end{displaymath}
After calculations, this equation reduces to
\begin{equation}\label{ce}
(\lambda+\delta)\left[\lambda+\delta+\beta(\overline{S})-(2\beta(\overline{S})+\alpha(\overline{N},\overline{S}))e^{-\delta\tau}e^{-\lambda\tau}\right]=0.
\end{equation}
We recall that the steady state $(\overline{S},\overline{N})$ is locally asymptotically stable when
all roots of (\ref{ce}) have negative real parts and the stability can only be lost if eigenvalues
cross the imaginary axis, that is if pure imaginary roots appear.

One can notice that $\lambda=-\delta<0$ is always an eigenvalue of (\ref{ce}). Therefore, we only
focus on the equation
\begin{equation}\label{ce2}
\lambda+\delta+\beta(\overline{S})-(2\beta(\overline{S})+\alpha(\overline{N},\overline{S}))e^{-\delta\tau}e^{-\lambda\tau}=0.
\end{equation}

We first analyze, in the next section, the stability of the trivial steady state $E^0$. We
establish necessary and sufficient conditions for the population's dying out. Then, in section
\ref{spss}, we concentrate on the behavior of the positive steady state $E^*$.

\section{Global Asymptotic Stability of the Trivial Steady State: Cell's Dying Out}\label{stss}

We concentrate, in this section, on the stability of the steady state $E^0=(0,0)$. From
(\ref{alpha}), $\alpha(0,0)=0$, so, for $\overline{S}=\overline{N}=0$, the characteristic equation
(\ref{ce2}) becomes
\begin{equation}\label{cet}
\lambda+\delta+\beta(0)-2\beta(0)e^{-\delta\tau}e^{-\lambda\tau}=0.
\end{equation}
It is straightforward to see that equation (\ref{cet}) has a unique real eigenvalue, say
$\lambda_0$, and all other eigenvalues $\lambda\neq \lambda_0$ of (\ref{cet}) satisfy
$\textrm{Re}(\lambda)<\lambda_0$.

Let consider the mapping $\lambda\mapsto
\lambda+\delta+\beta(0)-2\beta(0)e^{-\delta\tau}e^{-\lambda\tau}$ as a function of real $\lambda$.
Then it is an increasing function from $-\infty$ to $+\infty$, yielding the existence and
uniqueness of $\lambda_0$.

Assume that $\lambda=\mu+i\omega\neq\lambda_0$ satisfies (\ref{cet}). Then, considering the real
part of (\ref{cet}), we get
\begin{displaymath}
\mu-\lambda_0=2\beta(0)e^{-\delta\tau}\left[e^{-\mu\tau}\cos(\omega\tau)-e^{-\lambda_0\tau}\right].
\end{displaymath}
By contradiction, we assume that $\mu>\lambda_0$. Then
$e^{-\mu\tau}\cos(\omega\tau)-e^{-\lambda_0\tau}<0$ and we obtain a contradiction. So
$\mu\leq\lambda_0$. Now if $\mu=\lambda_0$, then the previous equality implies that
\begin{displaymath}
\cos(\omega\tau)=1, \qquad \textrm{ for } \tau\geq0.
\end{displaymath}
It follows that $\sin(\omega\tau)=0$ and, considering the imaginary part of (\ref{cet}) with
$\lambda=\mu+i\omega$, given by
\begin{displaymath}
\omega+2\beta(0)e^{-\delta\tau}e^{-\mu\tau}\sin(\omega\tau)=0,
\end{displaymath}
we obtain $\omega=0$ and $\lambda=\lambda_0$, which gives a contradiction. Therefore,
$\mu<\lambda_0$.

The real root $\lambda_0$ is negative if
\begin{displaymath}
(2e^{-\delta\tau}-1)\beta(0)<\delta,
\end{displaymath}
and all eigenvalues of (\ref{cet}) have negative real parts in this case. When condition
(\ref{condexist}) holds, $\lambda_0$ is positive. We can then conclude in the next proposition to
the stability of $E^0$.

\begin{proposition}\label{propzero}
The trivial steady state $E^0=(0,0)$ of system (\ref{eqS})--(\ref{eqN2}) is locally asymptotically
stable when
\begin{equation}\label{tss}
(2e^{-\delta\tau}-1)\beta(0)<\delta,
\end{equation}
and unstable when
\begin{displaymath}
(2e^{-\delta\tau}-1)\beta(0)>\delta.
\end{displaymath}
\end{proposition}

\begin{remark}
When
\begin{displaymath}
(2e^{-\delta\tau}-1)\beta(0)=\delta,
\end{displaymath}
then the unique real root of (\ref{cet}) is $\lambda_0=0$, and all other eigenvalues have negative
real parts. One can easily check that $\lambda_0=0$ is a simple root of (\ref{cet}), since the
first derivative of the mapping $\lambda\mapsto
\lambda+\delta+\beta(0)-2\beta(0)e^{-\delta\tau}e^{-\lambda\tau}$ at $\lambda_0=0$ is
\begin{displaymath}
1+2\beta(0)\tau e^{-\delta\tau}>0.
\end{displaymath}
Then the linear system is stable, but we cannot conclude to the asymptotic stability of the trivial
steady state $E^0=(0,0)$ of system (\ref{eqS})--(\ref{eqN2}) without further analysis. This is done
in Proposition \ref{propzero2}.
\end{remark}

When condition (\ref{tss}) holds true, $E^0$ is in fact the only steady state of system
(\ref{eqS})--(\ref{eqN2}) (see Proposition \ref{prop1}). In this case, we can show that $E^0$ is
globally asymptotically stable.

We first show the following result.

\begin{lemma}\label{lemproof}
Let $(S(t),N(t))$ be a solution of (\ref{eqS})--(\ref{eqN2}). If $\lim_{t\to+\infty}N(t)=0$, then
$\lim_{t\to +\infty}S(t)=0$.
\end{lemma}

\begin{proof}
Using (\ref{eqS}), a classical variation of constant formula gives us, for $t\geq0$,
\begin{displaymath}
S(t)=e^{-\delta t}S(0)+e^{-\delta t}\int_0^t
e^{\delta(\theta-\tau)}\beta(S(\theta-\tau))N(\theta-\tau) d\theta.
\end{displaymath}
Setting $\sigma=\theta-\tau$, this expression becomes
\begin{displaymath}
S(t)=e^{-\delta t}S(0) +e^{-\delta t}\int_{-\tau}^{t-\tau}
e^{\delta\sigma}\beta(S(\sigma))N(\sigma)d\sigma.
\end{displaymath}
Let $\varepsilon>0$ be fixed. Since $N$ is assumed to tend to zero when $t$ tends to $\infty$,
there exists $T>0$ such that
\begin{equation}\label{nbound}
N(t)<\varepsilon\frac{\delta e^{\delta\tau}}{2\beta(0)}, \qquad \textrm{ for } t\geq T.
\end{equation}
Then, for $t\geq T+\tau$,
\begin{displaymath}
S(t)=e^{-\delta t}\left[S(0)+\int_{-\tau}^{T} e^{\delta\sigma}\beta(S(\sigma))N(\sigma)d\sigma
\right] +e^{-\delta t}\int_{T}^{t-\tau} e^{\delta\sigma}\beta(S(\sigma))N(\sigma)d\sigma.
\end{displaymath}
Using (\ref{nbound}) and the fact that $\beta(0)$ is a bound of $\beta$, we obtain, for $t\geq
T+\tau$,
\begin{displaymath}
\begin{array}{rcl}
S(t) &\leq& e^{-\delta t}\left[S(0)+\displaystyle\int_{-\tau}^{T}
e^{\delta\sigma}\beta(S(\sigma))N(\sigma)d\sigma \right] +\varepsilon\displaystyle\frac{\delta
e^{\delta\tau}}{2}e^{-\delta t}\int_{T}^{t-\tau}
e^{\delta\sigma}d\sigma,\vspace{1ex}\\
&\leq& e^{-\delta t}\left[S(0)+\displaystyle\int_{-\tau}^{T}
e^{\delta\sigma}\beta(S(\sigma))N(\sigma)d\sigma \right] +\varepsilon
\displaystyle\frac{e^{\delta\tau}}{2}\left(
e^{-\delta\tau}-e^{-\delta(t-T)}\right),\vspace{1ex}\\
&\leq& e^{-\delta t}\left[S(0)+\displaystyle\int_{-\tau}^{T}
e^{\delta\sigma}\beta(S(\sigma))N(\sigma)d\sigma \right] +\displaystyle\frac{ \varepsilon}{2}.
\end{array}
\end{displaymath}
Let $\overline{t}>0$ be such that
\begin{displaymath}
e^{-\delta t}\left[S(0)+\displaystyle\int_{-\tau}^{T}
e^{\delta\sigma}\beta(S(\sigma))N(\sigma)d\sigma \right]<\frac{\varepsilon}{2}, \qquad \textrm{ for
} t\geq \overline{t}.
\end{displaymath}
Then, for $t\geq \max\{\overline{t}, T+\tau\}$, we obtain
\begin{displaymath}
S(t) < \varepsilon.
\end{displaymath}
Thus, $S(t)$ tends to zero as $t$ tends to $+\infty$, and the proof is complete.
\end{proof}

We recall a very useful lemma, proved by Barb\u alat (see Gopalsamy \cite{gopalsamy}).

\begin{lemma}\label{lemgop}
Let $f:[a,+\infty)\to\mathbb{R}$, $a\in\mathbb{R}$, be a differentiable function. If
$\lim_{t\to+\infty} f(t)$ exists and $f^{\prime}(t)$ is uniformly continuous on $(a,+\infty)$, then
\begin{displaymath}
\lim_{t\to+\infty} f^{\prime}(t)=0.
\end{displaymath}
\end{lemma}

We then prove the following result, dealing with the global asymptotic stability of $E^0$.

\begin{proposition}\label{propzero2}
Assume that
\begin{equation}\label{tssprop}
(2e^{-\delta\tau}-1)\beta(0)\leq \delta.
\end{equation}
Then all solutions $(S(t),N(t))$ of system (\ref{eqS})--(\ref{eqN2}) converge to the trivial
solution $(0,0)$. Hence $E^0$ is globally asymptotically stable and the cell populations dye out.
\end{proposition}

\begin{proof}
Let $(S(t),N(t))$ be a solution of (\ref{eqS})--(\ref{eqN2}). We define, for $t\geq0$,
\begin{displaymath}
Y(t)=N(t)+2e^{-\delta\tau}\int_{t-\tau}^t \beta(S(\theta))N(\theta) d\theta.
\end{displaymath}
Using (\ref{eqN2}), we can check that
\begin{equation}\label{Yp}
Y^{\prime}(t)=N(t)\left[(2e^{-\delta\tau}-1)\beta(S(t))-\delta\right].
\end{equation}
From condition (\ref{tssprop}) and the fact that $\beta$ is decreasing, it follows that
\begin{displaymath}
Y^{\prime}(t)\leq 0 \qquad \textrm{for } t>0.
\end{displaymath}
Thus, $Y$ is decreasing. Since $Y$ is a nonnegative function, we deduce that there exists $y\geq0$
such that
\begin{displaymath}
\lim_{t\to+\infty} Y(t)=y.
\end{displaymath}
In particularly, $Y$ is bounded, and consequently $N$ is also bounded.

We then deduce, with (\ref{eqN2}), that $N^{\prime}$ is bounded and, using a similar technic than
the one used in the proof of Lemma \ref{lemproof}, that $S$ is bounded. Consequently, with
(\ref{eqS}), we obtain that $S^{\prime}$ is bounded.

From (\ref{Yp}), since $N$, $S$, $N^{\prime}$ and $S^{\prime}$ are bounded, $Y^{\prime}$ is
uniformly continuous.

Since $\lim_{t\to+\infty} Y(t)$ exists and $Y^{\prime}$ is uniformly continuous on $(0,+\infty)$,
Lemma \ref{lemgop} implies that
\begin{displaymath}
\lim_{t\to+\infty} Y^{\prime}(t) = 0.
\end{displaymath}
Consequently, from (\ref{Yp}), we obtain either
\begin{displaymath}
\lim_{t\to+\infty} N(t)=0 \qquad \textrm{ or } \qquad \lim_{t\to+\infty}
(2e^{-\delta\tau}-1)\beta(S(t)) = \delta.
\end{displaymath}

First, assume that (\ref{tss}) holds true, that is
\begin{displaymath}
(2e^{-\delta\tau}-1)\beta(0)<\delta.
\end{displaymath}
If $2e^{-\delta\tau}-1>0$, then, since $\beta$ is a decreasing function satisfying (\ref{tss}), we
deduce that $(2e^{-\delta\tau}-1)\beta(S(t))\leq(2e^{-\delta\tau}-1)\beta(0)<\delta$. If
$2e^{-\delta\tau}-1\leq0$, then $(2e^{-\delta\tau}-1)\beta(S(t))\leq0<\delta$. Consequently, if it
exists, $\lim_{t\to+\infty} (2e^{-\delta\tau}-1)\beta(S(t))$ cannot be equal to $\delta$ and it
follows that
\begin{displaymath}
\lim_{t\to+\infty} N(t)=0.
\end{displaymath}
From Lemma \ref{lemproof}, we deduce that $\lim_{t\to+\infty}S(t)=0$, and the conclusion follows.

Second, assume that
\begin{displaymath}
(2e^{-\delta\tau}-1)\beta(0)=\delta.
\end{displaymath}
Then, $\lim_{t\to+\infty} (2e^{-\delta\tau}-1)\beta(S(t)) = \delta$ is equivalent to
$\lim_{t\to+\infty}\beta(S(t))=\beta(0)$. Since $\beta$ is positive and decreasing, this is
equivalent to $\lim_{t\to+\infty}S(t)=0$. It follows that either
\begin{displaymath}
\lim_{t\to+\infty} N(t)=0 \qquad \textrm{ or } \qquad \lim_{t\to+\infty} S(t)=0.
\end{displaymath}
If $\lim_{t\to+\infty} N(t)=0$, we conclude similarly to the previous case with Lemma
\ref{lemproof}. So we assume that $\lim_{t\to+\infty} S(t)=0$.

From (\ref{eqS}), we deduce that
\begin{displaymath}
\lim_{t\to+\infty}\beta(S(t-\tau))N(t-\tau)=0.
\end{displaymath}
Consequently, either
\begin{displaymath}
\lim_{t\to+\infty}\beta(S(t-\tau))=0 \qquad \textrm{ or } \qquad \lim_{t\to+\infty}N(t-\tau)=0.
\end{displaymath}
Since $\lim_{t\to+\infty} S(t)=0$, then $\lim_{t\to+\infty}\beta(S(t-\tau))=\beta(0)>0$. Hence,
$\lim_{t\to+\infty}N(t-\tau)=0$, and it follows that $\lim_{t\to+\infty}N(t)=0$.

This concludes the proof.
\end{proof}

\begin{remark}
Let $C$ denote the set of continuous functions mapping $[-\tau,0]$ into $\mathbb{R}^+$. One can
check that the function $V$, defined for $(\varphi,\psi)\in C^2$ by
\begin{displaymath}
V(\varphi,\psi)=\psi(0)+2e^{-\delta\tau}\int_{-\tau}^0 \beta(\varphi(\theta))\psi(\theta) d\theta,
\end{displaymath}
satisfies
\begin{displaymath}
\dot{V}(\varphi,\psi)=\psi(0)\left[(2e^{-\delta\tau}-1)\beta(\varphi(0))-\delta\right].
\end{displaymath}
Hence, $V$ is a Lyapunov functional (see Hale and Verduyn Lunel \cite{haleverduynlunel}) on the set
\begin{displaymath}
G=\left\{ (\varphi,\psi)\in C^2 ; \
\psi(0)\left[(2e^{-\delta\tau}-1)\beta(\varphi(0))-\delta\right]\leq 0 \right\}.
\end{displaymath}
With assumption (\ref{tssprop}), $G=C^2$. In the proof of Proposition \ref{propzero2}, we did not
directly use the properties of Lyapunov functionals, but the function $Y$ is defined by
\begin{displaymath}
Y(t)=V(S_t,N_t), \qquad \textrm{ for } t\geq0,
\end{displaymath}
where $S_t$ (respectively, $N_t$) is defined by $S_t(\theta)=S(t+\theta)$ (respectively,
$N_t(\theta)=N(t+\theta)$), $\theta\in[-\tau,0]$.
\end{remark}

Through Propositions \ref{propzero} and \ref{propzero2}, we obtained necessary and sufficient
conditions for the global asymptotic stability of $E^0$. Therefore, in the next section, we
concentrate on the behavior of $E^*$, the unique positive steady state of
(\ref{eqS})--(\ref{eqN2}).

\section{Local Asymptotic Stability of the Positive Steady State}\label{spss}

We now turn our considerations on the stability of the unique nontrivial steady state of system
(\ref{eqS})--(\ref{eqN2}), namely $E^*=(S^*,N^*)$, where, from Proposition \ref{prop1}, $S^*$ is
the unique solution of (\ref{eqss}), and $N^*=(2e^{-\delta\tau}-1)e^{\delta\tau}S^*$.

In order to ensure the existence of this steady state, we assume that condition (\ref{condexist}),
or equivalently condition (\ref{condexist2}), holds true. That is
\begin{displaymath}
\delta<\beta(0) \qquad \textrm{and} \qquad
0\leq\tau<\overline{\tau}:=\frac{1}{\delta}\ln\left(\frac{2\beta(0)}{\delta+\beta(0)}\right).
\end{displaymath}
In particularly, $2e^{-\delta\tau}-1>0$.

In this case, Proposition \ref{propzero} indicates that the unique other steady state $E^0$ is
unstable.

From their definitions in Proposition \ref{prop1}, the steady states $S^*$ and $N^*$ depend on the
time delay $\tau$. In fact,
\begin{displaymath}
S^*=S^*(\tau)=\beta^{-1}\left(\frac{\delta}{2e^{-\delta\tau}-1}\right) \quad \textrm{and}\quad
N^*=N^*(\tau)=\frac{2e^{-\delta\tau}-1}{e^{-\delta\tau}}S^*(\tau),
\end{displaymath}
where $\beta^{-1}: (0,\beta(0)]\to[0,+\infty)$ is a decreasing function.

Using these expressions, we can stress that $S^*$ and $N^*$ are positive decreasing continuous
functions of $\tau\in[0,\overline{\tau})$, continuously differentiable, such that
$S^*(0)=N^*(0)=\beta^{-1}(\delta)$, and
$\lim_{\tau\to\overline{\tau}}(S^*(\tau),N^*(\tau))=(0,0)=E^0$.

The characteristic equation (\ref{ce2}), with $\overline{S}=S^*$ and $\overline{N}=N^*$, is then
given by
\begin{equation}\label{cestar}
\lambda+A(\tau)-B(\tau)e^{-\lambda\tau}=0,
\end{equation}
with
\begin{displaymath}
A(\tau) := \delta+\beta(S^*(\tau)) \qquad \textrm{and} \qquad B(\tau) :=
[2\beta(S^*(\tau))+N^*(\tau)\beta^{\prime}(S^*(\tau))]e^{-\delta\tau}.
\end{displaymath}
Notice that $A(\tau)>0$ for all $\tau\in[0,\overline{\tau})$. Moreover, from (\ref{eqss}), we
obtain
\begin{equation}\label{B}
B(\tau)=A(\tau)+(2e^{-\delta\tau}-1)S^*(\tau)\beta^{\prime}(S^*(\tau)), \quad \textrm{for }
\tau\in[0,\overline{\tau}).
\end{equation}
In particular, $B(\tau)<A(\tau)$ for $\tau\in[0,\overline{\tau})$.

Taking $\tau=0$ in (\ref{cestar}), we obtain
\begin{displaymath}
\lambda+A(0)-B(0)=0,
\end{displaymath}
that is
\begin{displaymath}
\lambda=\beta^{-1}(\delta)\beta^{\prime}(\beta^{-1}(\delta)).
\end{displaymath}
Since $\beta$ is decreasing, we deduce that the only eigenvalue of (\ref{cestar}) is then negative.
The following lemma follows.

\begin{lemma}\label{lemmazero}
When $\delta<\beta(0)$ and $\tau=0$, the nontrivial steady-state $E^*$ of system
(\ref{eqS})--(\ref{eqN2}) is locally asymptotically stable, and the system
(\ref{eqS})--(\ref{eqN2}) undergoes a transcritical bifurcation.
\end{lemma}

When $\tau$ increases and remains in the interval $[0,\overline{\tau})$, the stability of the
steady state can only be lost if purely imaginary roots appear. Therefore, we investigate the
existence of purely imaginary roots of (\ref{cestar}).

Let $\lambda=i\omega$, $\omega\in\mathbb{R}$, be a pure imaginary eigenvalue of (\ref{cestar}).
Separating real and imaginary parts, we obtain
\begin{eqnarray}
A(\tau)-B(\tau)\cos(\omega\tau)&=&0,\label{eqcos}\\
\omega +B(\tau)\sin(\omega\tau)&=&0.\label{eqsin}
\end{eqnarray}

One can notice, firstly, that if $\omega$ is a solution of (\ref{eqcos})--(\ref{eqsin}) then
$-\omega$ also satisfies this system. Secondly, $\omega=0$ is not a solution of
(\ref{eqcos})--(\ref{eqsin}). Otherwise, we would obtain $A(\tau)=B(\tau)$ for some
$\tau\in[0,\overline{\tau})$, which contradicts $B(\tau)<A(\tau)$ for $\tau\in[0,\overline{\tau})$.
Therefore $\omega=0$ cannot be a solution of (\ref{eqcos})--(\ref{eqsin}). Thus, in the following,
we will only look for positive solutions $\omega$ of (\ref{eqcos})--(\ref{eqsin}).

From (\ref{eqcos}), a necessary condition for equation (\ref{cestar}) to have purely imaginary
roots is that
\begin{displaymath}
A(\tau)<|B(\tau)|.
\end{displaymath}
Since $B(\tau)<A(\tau)$, this implies in particularly that $B(\tau)$ must be negative. Moreover,
from (\ref{B}), the above condition is equivalent to
\begin{displaymath}
2A(\tau)+(2e^{-\delta\tau}-1)S^*(\tau)\beta^{\prime}(S^*(\tau))<0.
\end{displaymath}
Using the definitions of $A(\tau)$, $N^*(\tau)$, $S^*(\tau)$ and equality (\ref{ss1}), this
inequality becomes the following condition on $\tau$,
\begin{equation}\label{condtau}
\frac{4\delta e^{-\delta\tau}}{2e^{-\delta\tau}-1} +
(2e^{-\delta\tau}-1)\beta^{-1}\left(\frac{\delta}{2e^{-\delta\tau}-1}\right)\beta^{\prime}\left(\beta^{-1}\left(\frac{\delta}{2e^{-\delta\tau}-1}\right)\right)<0.
\end{equation}

\begin{lemma}\label{lemmachi}
Let $\chi : [0,+\infty)\to(-\infty,0]$ be defined, for $y\geq0$, by
\begin{displaymath}
\chi(y)=y\beta^{\prime}(y).
\end{displaymath}
Assume that
\begin{displaymath}
\begin{array}{ll}
\textrm{(H$_1$)} & \chi \textrm{ is decreasing on the interval } \left[0,\beta^{-1}(\delta)\right].
\vspace{1ex}\\
\textrm{(H$_2$)} & \chi\left(\beta^{-1}(\delta)\right)<-4\delta.
\end{array}
\end{displaymath}
Then there exists a unique $\tau^*\in(0,\overline{\tau})$ such that condition (\ref{condtau}) is
satisfied if and only if $\tau\in[0,\tau^*)$.
\end{lemma}

\begin{proof}
Consider the negative functions $f_1(\tau)$ and $f_2(\tau)$, defined for
$\tau\in[0,\overline{\tau}]$ by
\begin{displaymath}
f_1(\tau)=\chi\left(\beta^{-1}\left(\frac{\delta}{2e^{-\delta\tau}-1}\right)\right) \qquad
\textrm{and} \qquad f_2(\tau)=-\frac{4\delta e^{-\delta\tau}}{(2e^{-\delta\tau}-1)^2}.
\end{displaymath}
Then
\begin{displaymath}
\left\{\tau\in[0,\overline{\tau}) \ ; \textrm{condition } (\ref{condtau}) \textrm{ is satisfied}
\right\} = \left\{\tau\in[0,\overline{\tau}) \ ; f_1(\tau)<f_2(\tau) \right\}.
\end{displaymath}

The function $f_2$ satisfies
\begin{displaymath}
f_2(0)=-4\delta \qquad \textrm{ and } \qquad
f_2(\overline{\tau})=-2\frac{\delta+\beta(0)}{\delta}\beta(0),
\end{displaymath}
and, for $\tau\in[0,\overline{\tau}]$,
\begin{displaymath}
f_2^{\prime}(\tau)=-\frac{4\delta^2 e^{-\delta\tau}
(2e^{-\delta\tau}+1)}{(2e^{-\delta\tau}-1)^3}<0.
\end{displaymath}
Hence $f_2$ is decreasing from $-4\delta$ to $-2(\delta+\beta(0))\beta(0)/\delta$.

For $\tau\in[0,\overline{\tau})$,
\begin{displaymath}
\delta\leq\frac{\delta}{2e^{-\delta\tau}-1}<\beta(0).
\end{displaymath}
Since $\beta^{-1}$ is decreasing on $(0,\beta(0)]$ and, from (H$_1$), $\chi$ is decreasing on
$[0,\beta^{-1}(\delta)]$, we deduce that $f_1$ is increasing. Moreover,
\begin{displaymath}
f_1(0)=\chi\left(\beta^{-1}\left(\delta\right)\right) \qquad \textrm{ and } \qquad
f_1(\overline{\tau})=0.
\end{displaymath}
From (H$_2$),
\begin{displaymath}
f_1(0)<f_2(0).
\end{displaymath}
Since $f_1(\overline{\tau})>f_2(\overline{\tau})$, there exists $\tau^*\in(0,\overline{\tau})$,
which is unique since $f_1$ is increasing and $f_2$ decreasing, such that
\begin{displaymath}
\left\{\tau\in[0,\overline{\tau}) \ ; f_1(\tau)<f_2(\tau) \right\} = [0,\tau^*).
\end{displaymath}
This concludes the proof.
\end{proof}

\begin{remark}\label{remchi}
Assumption (\textit{H}$_1$) is not necessary for the existence of $\tau^*$, it just implies the
uniqueness. This latter can be achieved with weaker conditions, as we will check on an example in
section \ref{snum}.
\end{remark}

We assume, in the following, that (\textit{H}$_1$) and (\textit{H}$_2$) are fulfilled and
$\tau\in[0,\tau^*)$.

System (\ref{eqcos})--(\ref{eqsin}) is equivalent to
\begin{equation}
\cos(\omega\tau)=\frac{A(\tau)}{B(\tau)}, \qquad
\sin(\omega\tau)=-\frac{\omega}{B(\tau)}.\label{eqcossin}
\end{equation}
Note that for $\tau\in[0,\tau^*)$, $B(\tau)<0$.

Therefore, adding the squares of both sides of (\ref{eqcossin}), purely imaginary eigenvalues
$i\omega$ of (\ref{cestar}), with $\omega>0$, must satisfy
\begin{equation}\label{omega}
\omega=\sqrt{B^2(\tau)-A^2(\tau)}.
\end{equation}
So, in the following, we will think of $\omega$ as $\omega(\tau)$.

Substituting this expression for $\omega$ in (\ref{eqcossin}), we obtain
\begin{equation}\label{eqcossin2}
\begin{array}{rcl}
\cos\left(\tau\sqrt{B^2(\tau)-A^2(\tau)}\right)&=&\displaystyle\frac{A(\tau)}{B(\tau)}, \vspace{1ex}\\
\sin\left(\tau\sqrt{B^2(\tau)-A^2(\tau)}\right)&=&-\displaystyle\frac{\sqrt{B^2(\tau)-A^2(\tau)}}{B(\tau)}.
\end{array}
\end{equation}
From the above reasoning, values of $\tau\in[0,\tau^*)$ solutions of system (\ref{eqcossin2})
generate positive $\omega(\tau)$, given by (\ref{omega}), and hence yield imaginary eigenvalues of
(\ref{cestar}). Consequently, we look for positive solutions $\tau$ of (\ref{eqcossin2}) in the
interval $(0,\tau^*)$.

Positive solutions $\tau\in(0,\tau^*)$ of (\ref{eqcossin2}) satisfy
\begin{displaymath}
\tau\sqrt{B^2(\tau)-A^2(\tau)} = \arccos\left(\frac{A(\tau)}{B(\tau)}\right) + 2k\pi, \qquad
k\in\mathbb{N}_0,
\end{displaymath}
where $\mathbb{N}_0$ denotes the set of all nonnegative integers. We set
\begin{displaymath}
\tau_k(\tau)=\frac{\arccos\left(\frac{A(\tau)}{B(\tau)}\right) +
2k\pi}{\sqrt{B^2(\tau)-A^2(\tau)}}, \qquad k\in\mathbb{N}_0, \tau\in[0,\tau^*).
\end{displaymath}
Values of $\tau$ for which $\omega(\tau)=\sqrt{B^2(\tau)-A^2(\tau)}$ is a solution of
(\ref{eqcossin}) are roots of the functions
\begin{equation}\label{Sn}
Z_k(\tau)=\tau-\tau_k(\tau), \qquad k\in\mathbb{N}_0, \tau\in[0,\tau^*).
\end{equation}
The roots of $Z_k$ can be found using popular software like Maple, but are hard to determine with
analytical tools \cite{bk2002}. The following lemma states some properties of the $Z_k$ functions.

\begin{lemma}\label{lemmaZ}
For $k\in\mathbb{N}_0$,
\begin{displaymath}
Z_k(0)<0 \qquad \textrm{ and } \qquad \lim_{\tau\to\tau^*}Z_k(\tau)=-\infty.
\end{displaymath}
Therefore, provided that no root of $Z_k$ is a local extremum, the number of positive roots of
$Z_k$, $k\in\mathbb{N}_0$, on the interval $[0,\tau^*)$ is even.

Moreover, if $Z_k$ has no root on the interval $[0,\tau^*)$, then $Z_j$, with $j>k$, does not have
positive roots.
\end{lemma}

\begin{proof}
Notice first that $\tau_k(0)>0$ and, secondly, that
\begin{displaymath}
\lim_{\tau\to\tau^*}\arccos\left(\frac{A(\tau)}{B(\tau)}\right)=\pi \qquad\textrm{and}\qquad
\lim_{\tau\to\tau^*}\sqrt{B^2(\tau)-A^2(\tau)}=0,
\end{displaymath}
so $\lim_{\tau\to\tau^*}\tau_k(\tau)=+\infty$. Then the first statement holds.

To prove the second statement, one can notice that $\tau_{k+1}(\tau)>\tau_{k}(\tau)$, for
$\tau\in[0,\tau^*)$ and $k\in\mathbb{N}_0$, so
\begin{displaymath}
Z_{k+1}(\tau)<Z_k(\tau), \qquad \tau\in[0,\tau^*), k\in\mathbb{N}_0.
\end{displaymath}
Using the fact that $Z_k(0)<0$, we conclude. This ends the proof.
\end{proof}

\begin{remark}\label{rem1}
The second statement in Lemma \ref{lemmaZ} implies, in particularly, that, if $Z_0$ has no positive
root, then (\ref{eqcossin}) has no positive solution, and equation (\ref{cestar}) does not have
pure imaginary roots.
\end{remark}

In the following proposition we establish some properties of pure imaginary roots of equation
(\ref{cestar}), using a method described in \cite{kuang1993}.

\begin{proposition}\label{propomega}
Let $\pm i\omega(\tau_c)$, with $\omega(\tau_c)>0$, be a pair of pure imaginary roots of equation
(\ref{cestar}) when $\tau=\tau_c$. Then $\pm i\omega(\tau_c)$ are simple roots of (\ref{cestar})
such that
\begin{equation}\label{sign}
\textrm{sign}\left\{\frac{d\textrm{Re}(\lambda)}{d\tau}\bigg|_{\tau=\tau_c}\right\} =
\textrm{sign}\bigg\{
-B^3-B^2B^{\prime}\tau_c+B(A^2+A^{\prime}+AA^{\prime}\tau_c)-B^{\prime}A\bigg\},
\end{equation}
where $A=A(\tau_c)$, $B=B(\tau_c)$, $A^{\prime}=A^{\prime}(\tau_c)$ and
$B^{\prime}=B^{\prime}(\tau_c)$.
\end{proposition}

\begin{proof}
We set
\begin{displaymath}
\Delta(\lambda,\tau)=\lambda+A(\tau)-B(\tau)e^{-\lambda\tau}.
\end{displaymath}
Let $\lambda(\tau)$ be a family of roots of (\ref{cestar}), so $\Delta(\lambda(\tau),\tau)=0$, such
that $\lambda(\tau_c)$ is a pure imaginary root of (\ref{cestar}), given by
$\lambda(\tau_c)=i\omega(\tau_c)$. Then,
\begin{equation}\label{eqdelta}
\frac{d\lambda}{d\tau}(\tau)\Delta_{\lambda}(\lambda,\tau) + \Delta_{\tau}(\lambda,\tau) = 0,
\end{equation}
where
\begin{displaymath}
\Delta_{\lambda}(\lambda,\tau) := \frac{d\Delta}{d\lambda}(\lambda,\tau) = 1+B(\tau)\tau
e^{-\lambda\tau},
\end{displaymath}
and
\begin{displaymath}
\Delta_{\tau}(\lambda,\tau) := \frac{d\Delta}{d\tau}(\lambda,\tau) =
A^{\prime}(\tau)+[B(\tau)\lambda -B^{\prime}(\tau)]e^{-\lambda\tau}.
\end{displaymath}

Assume, by contradiction, that $\lambda(\tau_c)=i\omega(\tau_c)$ is not a simple root of
(\ref{cestar}). Then, from (\ref{eqdelta}),
$\Delta_{\tau}(i\omega(\tau_c),\tau_c)=A^{\prime}(\tau_c)+[iB(\tau_c)\omega(\tau_c)
-B^{\prime}(\tau_c)]e^{-i\omega(\tau_c)\tau_c}=0$. Separating real and imaginary parts in this
equality we deduce
\begin{displaymath}
\begin{array}{rcl}
B^{\prime}(\tau_c)\cos(\omega(\tau_c)\tau_c)-B(\tau_c)\omega(\tau_c)\sin(\omega(\tau_c)\tau_c)&=&A^{\prime}(\tau_c),\vspace{1ex}\\
B(\tau_c)\omega(\tau_c)\cos(\omega(\tau_c)\tau_c)+B^{\prime}(\tau_c)\sin(\omega(\tau_c)\tau_c)&=&0.
\end{array}
\end{displaymath}
We recall that $B(\tau_c)$ is necessarily strictly negative. Using (\ref{eqcossin}), the above
system is equivalent to
\begin{displaymath}
\begin{array}{rcl}
\omega(\tau_c)^2&=&\displaystyle\frac{A(\tau_c)B^{\prime}(\tau_c)-B(\tau_c)A^{\prime}(\tau_c)}{B(\tau_c)},\vspace{1ex}\\
\omega(\tau_c)\left[A(\tau_c)+\displaystyle\frac{B^{\prime}(\tau_c)}{B(\tau_c)}\right]&=&0.
\end{array}
\end{displaymath}
Since $\omega(\tau_c)>0$ and satisfies, from (\ref{omega}),
$\omega(\tau_c)^2=B(\tau_c)^2-A(\tau_c)^2$, we obtain
\begin{displaymath}
\begin{array}{rcl}
B(\tau_c)^2-A(\tau_c)^2&=&\displaystyle\frac{A(\tau_c)B^{\prime}(\tau_c)-B(\tau_c)A^{\prime}(\tau_c)}{B(\tau_c)},\vspace{1ex}\\
A(\tau_c)B(\tau_c)&=&-B^{\prime}(\tau_c).
\end{array}
\end{displaymath}
Substituting the second equation in the first one, this yields
$B(\tau_c)^2-A(\tau_c)^2=-A^{\prime}(\tau_c)-A(\tau_c)^2$, so
\begin{displaymath}
B(\tau_c)^2=-A^{\prime}(\tau_c).
\end{displaymath}
Since $B(\tau_c)^2>0$ and $A^{\prime}(\tau_c)>0$, we obtain a contradiction. Hence
$\lambda(\tau_c)=i\omega(\tau_c)$ is a simple root of (\ref{cestar}).

In the following, we do not mention the dependence of the coefficients $A$ and $B$ (and their
derivatives) with respect to $\tau$.

Now, from (\ref{eqdelta}), we obtain
\begin{displaymath}
\left(\frac{d\lambda}{d\tau}\right)^{-1} = \frac{e^{\lambda\tau}+B\tau}{B^{\prime}-B\lambda
-A^{\prime}e^{\lambda\tau}}.
\end{displaymath}
Since $\Delta(\lambda,\tau)=0$, we deduce
\begin{displaymath}
e^{\lambda\tau} = \frac{B}{\lambda+A}.
\end{displaymath}
Therefore,
\begin{displaymath}
\left(\frac{d\lambda}{d\tau}\right)^{-1} =
\frac{B+B\tau(\lambda+A)}{(B^{\prime}-B\lambda)(\lambda+A) -A^{\prime}B}.
\end{displaymath}
For $\tau=\tau_c$, we obtain
\begin{displaymath}
\begin{array}{rcl}
\left(\displaystyle\frac{d\lambda}{d\tau}\right)^{-1}\bigg|_{\tau=\tau_c} &=&
\displaystyle\frac{B+B\tau_c(i\omega(\tau_c)+A)}{(B^{\prime}-iB\omega(\tau_c))(i\omega(\tau_c)+A)
-A^{\prime}B},\vspace{1ex}\\
&=&\displaystyle\frac{B(1+A\tau_c)+iB\omega(\tau_c)\tau_c}{B^{\prime}A-AB^{\prime}+B\omega^2(\tau_c)+i(B^{\prime}-AB)\omega(\tau_c)}.
\end{array}
\end{displaymath}
Then,
\begin{displaymath}
\textrm{Re}\left(\frac{d\lambda}{d\tau}\right)^{-1}\bigg|_{\tau=\tau_c} =
\frac{[B^2(1+A\tau_c)+B\tau_c(B^{\prime}-AB)]\omega(\tau_c)^2+B(1+A\tau_c)(B^{\prime}A-A^{\prime}B)
}{[B^{\prime}A-A^{\prime}B+B\omega(\tau_c)^2] + [B^{\prime}-AB]^2\omega(\tau_c)^2}.
\end{displaymath}

Noticing that
\begin{displaymath}
\textrm{sign}\left\{ \frac{d\textrm{Re}(\lambda)}{d\tau}\right\} = \textrm{sign}\left\{
\textrm{Re}\left(\frac{d\lambda}{d\tau}\right)^{-1}\right\},
\end{displaymath}
we get,
\begin{equation}\label{eqtemp}
\begin{array}{rcl}
\textrm{sign}\left\{ \displaystyle\frac{d\textrm{Re}(\lambda)}{d\tau}\bigg|_{\tau=\tau_c}\right\}
&=& \textrm{sign}\bigg\{[B^2(1+A\tau_c)+B\tau_c(B^{\prime}-AB)]\omega(\tau_c)^2 \\
&&\qquad\qquad\qquad\quad +B(1+A\tau_c)(B^{\prime}A-A^{\prime}B)\bigg\}.
\end{array}
\end{equation}
Since $i\omega(\tau_c)$ is a purely imaginary root of (\ref{cestar}), then, from (\ref{omega}),
\begin{displaymath}
\omega(\tau_c)^2=B^2-A^2.
\end{displaymath}
Substituting this expression in (\ref{eqtemp}), we obtain, after simplifications
\begin{displaymath}
\textrm{sign}\left\{ \displaystyle\frac{d\textrm{Re}(\lambda)}{d\tau}\bigg|_{\tau=\tau_c}\right\} =
\textrm{sign}\bigg\{B\left[B^3+B^2B^{\prime}\tau_c-B(A^2+A^{\prime}+AA^{\prime}\tau_c)+B^{\prime}A\right]\bigg\}.
\end{displaymath}
As we already noticed, if equation (\ref{cestar}) has pure imaginary roots then necessarily $B<0$.
We then deduce (\ref{sign}) and the proof is complete.
\end{proof}

Using this last proposition and the previous results about the existence of purely imaginary roots
of (\ref{ce}), we can state and prove the following theorem, dealing with the asymptotic stability
of $E^*$.

\begin{theorem}\label{theohopf}
Assume that (\ref{condexist2}) holds true and (H$_1$) and (H$_2$) are fulfilled.
\begin{itemize}
  \item[(i)] If $Z_0$ (defined in (\ref{Sn})) has no root on the interval $[0,\tau^*)$,
$\tau^*$ defined in Lemma \ref{lemmachi}, then the positive steady state $E^*=(S^*,N^*)$ of
(\ref{eqS})--(\ref{eqN2}) is locally asymptotically stable for $\tau\in[0,\overline{\tau})$.
  \item[(ii)] If $Z_0$ has at least one positive root $\tau_c\in(0,\tau^*)$ then $E^*$ is
locally asymptotically stable for $\tau\in[0,\tau_c)$ and a Hopf bifurcation occurs at $E^*$ for
$\tau=\tau_c$ if
\begin{equation}\label{tc}
B(A^2+A^{\prime}+AA^{\prime}\tau_c)-B^3-B^2B^{\prime}\tau_c-B^{\prime}A\neq 0,
\end{equation}
where $A=A(\tau_c)$, $B=B(\tau_c)$, $A^{\prime}=A^{\prime}(\tau_c)$ and
$B^{\prime}=B^{\prime}(\tau_c)$.
\end{itemize}
\end{theorem}

\begin{proof}
First, from Lemma \ref{lemmazero}, we know that $E^*$ is locally asymptotically stable when
$\tau=0$.

If $Z_0$ has no positive root on the interval $(0,\tau^*)$, then the characteristic equation
(\ref{ce}) has no pure imaginary root (see Remark \ref{rem1} and Lemma \ref{lemmaZ}). Consequently,
the stability of $E^*$ cannot be lost when $\tau$ increases. We obtain the statement in
\textbf{(i)}.

Now, if $Z_0$ has at least one positive root, say $\tau_c\in(0,\tau^*)$, then equation (\ref{ce})
has a pair of simple conjugate pure imaginary roots $\pm i\omega(\tau_c)$ for $\tau=\tau_c$. From
(\ref{tc}) together with Proposition \ref{propomega}, we have either
\begin{displaymath}
\displaystyle\frac{d\textrm{Re}(\lambda)}{d\tau}\bigg|_{\tau=\tau_c}>0 \qquad \textrm{ or } \qquad
\displaystyle\frac{d\textrm{Re}(\lambda)}{d\tau}\bigg|_{\tau=\tau_c}<0.
\end{displaymath}
By contradiction, we assume that there exists a branch of characteristic roots $\lambda(\tau)$ such
that $\lambda(\tau_c)=i\omega_c$ and
\begin{displaymath}
\frac{d\textrm{Re}(\lambda(\tau))}{d\tau}<0
\end{displaymath}
for $\tau<\tau_c$, $\tau$ close to $\tau_c$. Then there exists a characteristic root
$\lambda(\tau)$ such that $\textrm{Re}(\lambda(\tau))>0$ and $\tau<\tau_c$. Since $E^*$ is locally
asymptotically stable when $\tau=0$, applying Rouch\'e's Theorem \cite{dieudonne}, we obtain that
all characteristic roots of (\ref{ce}) have negative real parts when $\tau\in[0,\tau_c)$, and we
obtain a contradiction. Thus,
\begin{displaymath}
\displaystyle\frac{d\textrm{Re}(\lambda)}{d\tau}\bigg|_{\tau=\tau_c}>0.
\end{displaymath}
In this case, a Hopf bifurcation occurs at $E^*$ when $\tau=\tau_c$.
\end{proof}

The result stated in (ii) leads, through the Hopf bifurcation, to the existence of periodic
solutions for system (\ref{eqS})--(\ref{eqN2}).

In the next section, we apply the above-mentioned results of stability to a particular introduction
rate $\beta$ and we present some numerical illustrations.

\section{Example and Numerical Simulations}\label{snum}

We develop, in this section, numerical illustrations of the above mentioned results (mainly the
ones stated in Theorem \ref{theohopf}).

Let define (see \cite{m1978, m1979, pbm2005, pm2004}) the introduction  rate $\beta$ by
\begin{displaymath}
\beta(S)=\beta_0\frac{\theta^n}{\theta^n+S^n}, \qquad \beta_0, \theta \geq0, \ n>1.
\end{displaymath}
The parameter $\beta_0$ represents the maximal rate of introduction in the proliferating phase,
$\theta$ is the value for which $\beta$ attains half of its maximum value, and $n$ is the
sensitivity of the rate of reintroduction. The coefficient $n$ describes the reaction of $\beta$
due to external stimuli, the action of a growth factor for example (some growth factors are known
to trigger the introduction of nonproliferating cells in the proliferating phase \cite{bmm1995,
mbm1998}).

Then, from (\ref{condexist2}), the unique positive steady state of (\ref{eqS})--(\ref{eqN2}) exists
if and only if
\begin{displaymath}
\delta<\beta_0 \qquad \textrm{and} \qquad
0\leq\tau<\overline{\tau}=\frac{1}{\delta}\ln\left(\frac{2\beta_0}{\delta+\beta_0}\right).
\end{displaymath}
From (\ref{ss1})--(\ref{ss2}), it is defined by
\begin{displaymath}
S^*=\theta\left(\frac{(2e^{-\delta\tau}-1)\beta_0}{\delta}-1\right)^{1/n} \quad \textrm{and} \quad
N^*=\theta\frac{2e^{-\delta\tau}-1}{e^{-\delta\tau}}\left(\frac{(2e^{-\delta\tau}-1)\beta_0}{\delta}-1\right)^{1/n}
\end{displaymath}
Note that the function $\beta^{-1}$ is defined by $\beta^{-1}(x)=\theta(\beta_0/x-1)^{1/n}$ for
$x\in(0,\beta_0]$.

After computations, we can state that condition (\ref{condtau}) is equivalent to
\begin{displaymath}
\tau<\tau^* := \frac{1}{\delta}\ln\left(\frac{2\beta_0(n-2)}{n(\beta_0+\delta)}\right),
\end{displaymath}
provided that
\begin{equation}\label{n}
n>\frac{4\beta_0}{\beta_0-\delta}.
\end{equation}
Noticing that the function $\chi$, defined in Lemma \ref{lemmachi}, is given by
\begin{displaymath}
\chi(y)=-n\beta_0\frac{\theta^ny^n}{(\theta^n+y^n)^2},
\end{displaymath}
one can check that (\ref{n}) is equivalent to (\textit{H}$_2$).

However, the function $\chi$ does not necessarily satisfy (\textit{H}$_1$), which is too strong (as
mentioned in Remark \ref{remchi}). For $y\geq0$,
\begin{displaymath}
\chi^{\prime}(y)=\frac{\beta_0n^2\theta^ny^{n-1}}{(\theta^n+y^n)^3}(y^n-\theta^n).
\end{displaymath}
Consequently, $\chi$ is decreasing for $y\leq\theta$ and increasing for $y>\theta$. Taking
$y=\beta^{-1}(\delta)$, we find that $\chi$ is decreasing on $[0,\beta^{-1}(\delta)]$ if and only
if $\beta_0<2\delta$. In this case (\textit{H}$_1$) is fulfilled. If $\beta_0>2\delta$, then $\chi$
is decreasing on the interval $[0,\theta]$ and increasing on $[\theta,\beta^{-1}(\delta)]$, yet
$\tau^*$ is uniquely defined.

Note that
\begin{displaymath}
A(\tau)=\frac{2\delta e^{-\delta\tau}}{2e^{-\delta\tau}-1} \quad \textrm{ and } \quad
B(\tau)=\frac{2\delta\beta_0
e^{-\delta\tau}-n\delta[(2e^{-\delta\tau}-1)\beta_0-\delta]}{(2e^{-\delta\tau}-1)\beta_0},
\end{displaymath}
with
\begin{displaymath}
A^{\prime}(\tau)=\frac{2\delta^2 e^{-\delta\tau}}{(2e^{-\delta\tau}-1)^2} \qquad \textrm{ and }
\qquad B^{\prime}(\tau)=\frac{2\delta^2
e^{-\delta\tau}(\beta_0+n\delta)}{(2e^{-\delta\tau}-1)^2\beta_0}.
\end{displaymath}
Then, assuming that (\ref{n}) holds true, we define the $Z_k$ functions, as in (\ref{Sn}), for
$\tau\in[0,\tau^*)$ by
\begin{displaymath}
Z_k(\tau)=\tau-\frac{\arccos\left(\frac{A(\tau)}{B(\tau)}\right) +
2k\pi}{\sqrt{B^2(\tau)-A^2(\tau)}}, \quad k\in\mathbb{N}_0, \tau\in[0,\tau^*).
\end{displaymath}

We choose the parameters according to \cite{m1978, pbm2005, pm2004}:
\begin{equation}\label{par}
\delta=0.05 \textrm{ days}^{-1}, \quad \beta_0=1.77 \textrm{ days}^{-1}, \quad \theta=1.
\end{equation}
Notice that the value of $\theta$ is in fact normalized and does not influence the stability of
system (\ref{eqS})--(\ref{eqN2}) since all coefficients actually do not depend on $\theta$. The
value of $\theta$ only influences the shape of the oscillations and the values of the steady
states.

Using Maple to determine the roots of $Z_n$, we first check that $Z_0$ (and consequently all $Z_k$
functions) is strictly negative on $[0,\tau^*)$ for $n\leq 10$. Hence, from Theorem \ref{theohopf},
the positive steady state $E^*=(S^*,N^*)$ of (\ref{eqS})--(\ref{eqN2}) is locally asymptotically
stable for $\tau\in[0,\overline{\tau})$.

For $n\geq 10$, Pujo-Menjouet et al. \cite{pbm2005, pm2004} noticed, for the model
(\ref{eqP})--(\ref{eqN}) with the introduction rate $\beta$ depending only upon the
nonproliferating phase population $N(t)$, that oscillations may be observed.

\begin{figure}[pt]
\begin{center}
\includegraphics[width=6cm, height=4cm]{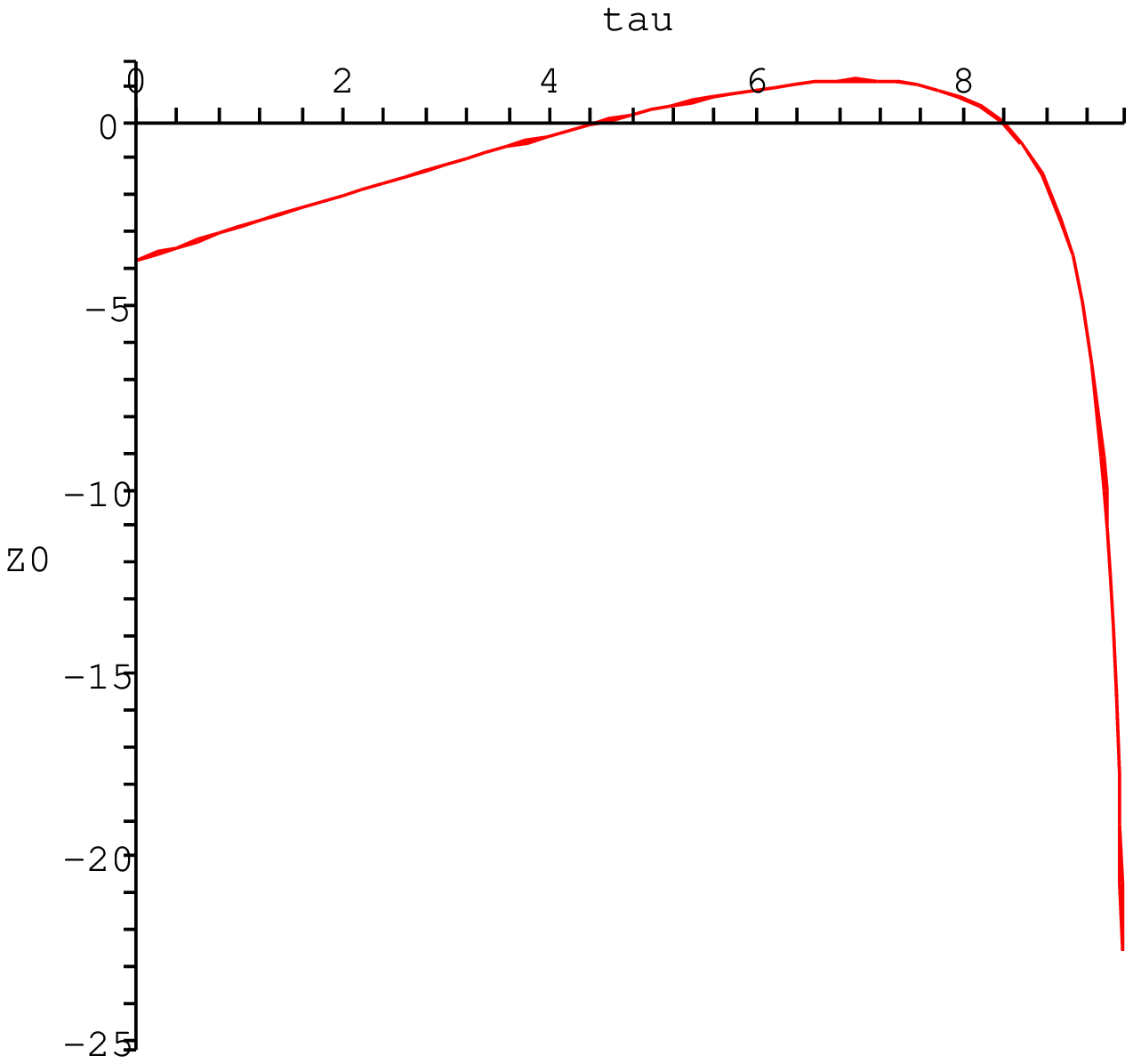}
\includegraphics[width=6cm, height=4cm]{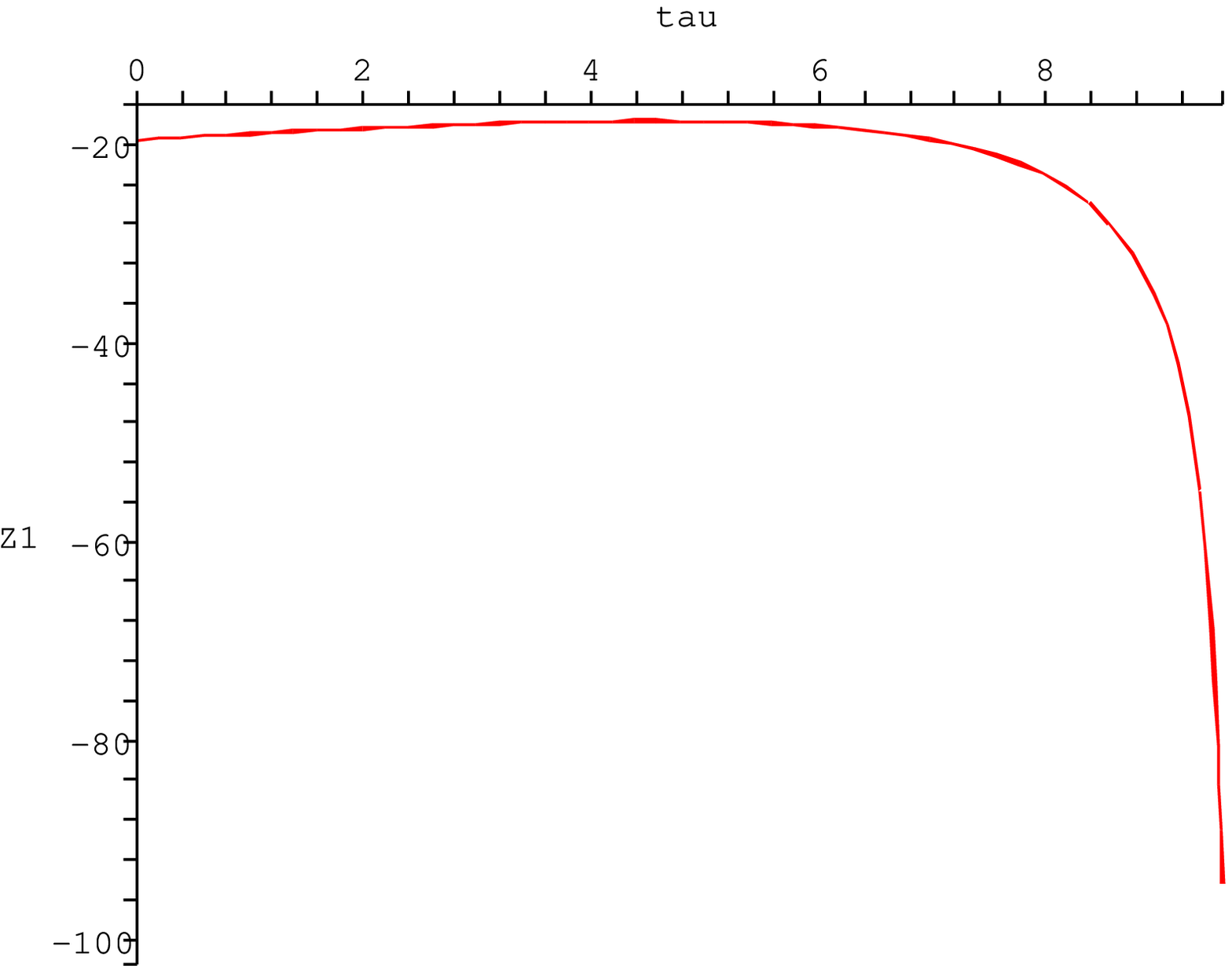}
\caption{The functions $Z_0$ (left) and $Z_1$ (right) are drawn on the interval $[0, \tau^*)$ for
parameters given by (\ref{par}) and $n=12$. One can see that $Z_0$ has exactly two roots,
$\tau_1\simeq 4.52$ and $\tau_2\simeq 8.36$, and $Z_1$ has no root.}\label{Z}
\end{center}
\end{figure}

We choose $n=12$, in keeping with values in \cite{pbm2005, pm2004}. Then, we find that
\begin{displaymath}
\overline{\tau}\simeq 13.3 \textrm{ days} \qquad \textrm{ and } \qquad \tau^*\simeq 9.66 \textrm{
days}.
\end{displaymath}
One can see on Figure \ref{Z} that $Z_0$ has two positive roots in this case, $\tau_1\simeq 4.52$
days and $\tau_2\simeq 8.36$ days, and that $Z_1$ is strictly negative, so all $Z_k$ functions,
with $k\geq1$ have no roots. Consequently, there exist two critical values, $\tau_1$ and $\tau_2$,
for which a stability switch can occur at $E^*$.

\begin{figure}[pt]
\begin{center}
\includegraphics[width=6cm, height=4cm]{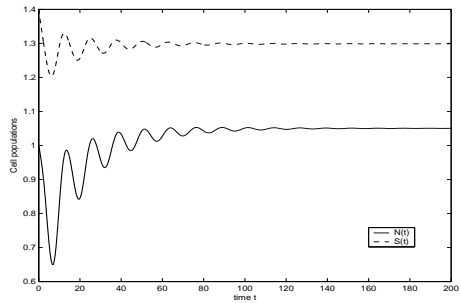}
\caption{For $\tau=3.5$ days, and the other parameters given by (\ref{par}) with $n=12$, the
solutions $S(t)$ (dashed line) and $N(t)$ (solid line) oscillate transiently to the steady state,
which is asymptotically stable. Damped oscillations are observed.}\label{stab3p5}
\end{center}
\end{figure}

For $\tau<\tau_1$, one can check that the populations are asymptotically stable on Figure
\ref{stab3p5}. In this case $\tau=3.5$ days and the solutions of (\ref{eqS})--(\ref{eqN2})
oscillate transiently to the steady state. Numerical simulations of the solutions of
(\ref{eqS})--(\ref{eqN2}) are carried out with dde23 \cite{dde23}, a {\sc Matlab} solver for delay
differential equations.

\begin{figure}[pt]
\begin{center}
\begin{subfigure}[]{
\includegraphics[width=6cm, height=4cm]{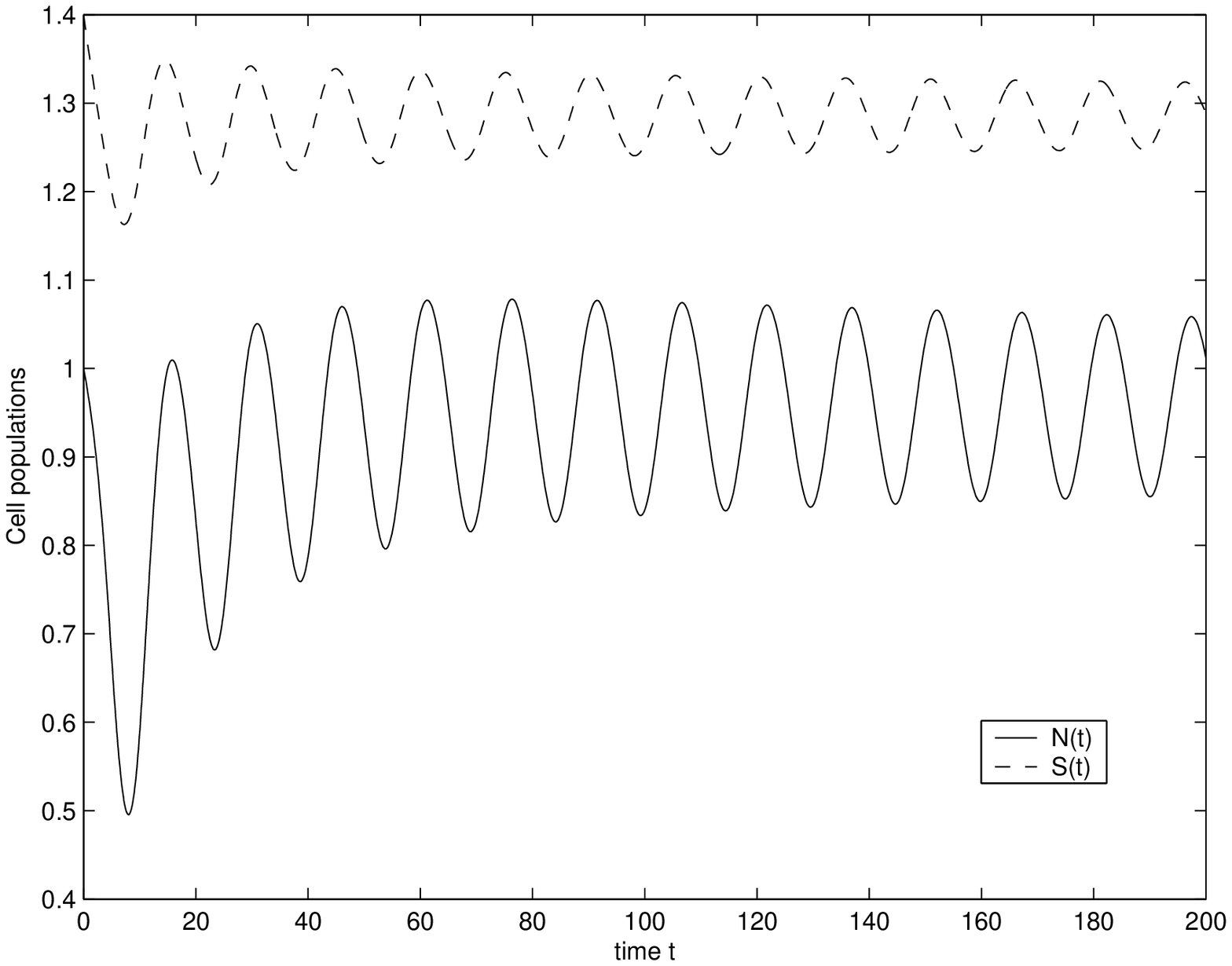}
}\end{subfigure}
\begin{subfigure}[]{
\includegraphics[width=6cm, height=4cm]{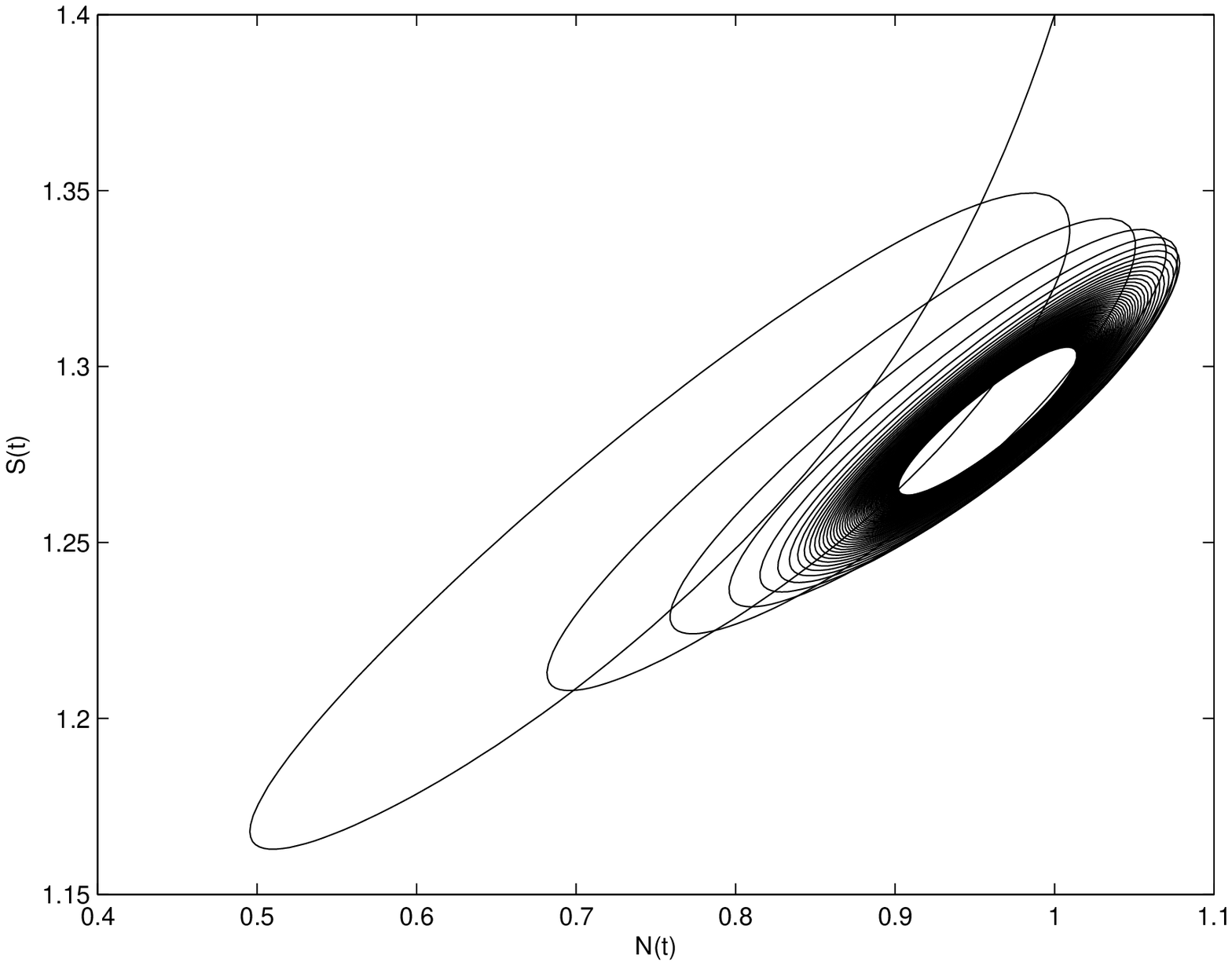}
}\end{subfigure} \caption{For $\tau=4.52$ days, and the other parameters given by (\ref{par}) with
$n=12$, a Hopf bifurcation occurs and the steady state $(S^*,N^*)$ of (\ref{eqS})--(\ref{eqN2}) is
unstable. The periodic solutions $S(t)$ (dashed line) and $N(t)$ (solid line) are represented in
(a), and we can observe the solutions in the $(S,N)$-plane in (b). Periods of the oscillations are
about 15 days. }\label{bif1}
\end{center}
\end{figure}

When $\tau=\tau_1$, one can check that
\begin{displaymath}
B(A^2+A^{\prime}+AA^{\prime}\tau_c)-B^3-B^2B^{\prime}\tau_c-B^{\prime}A\simeq 0.053,
\end{displaymath}
so condition (\ref{tc}) holds, and a Hopf bifurcation occurs at $(S^*,N^*)$, from Theorem
\ref{theohopf}. This is illustrated on Figure \ref{bif1}. Periodic solutions with periods about
$15$ days are observed at the bifurcation, and the steady state $E^*$ becomes unstable.

\begin{figure}[pt]
\begin{center}
\includegraphics[width=6cm, height=4cm]{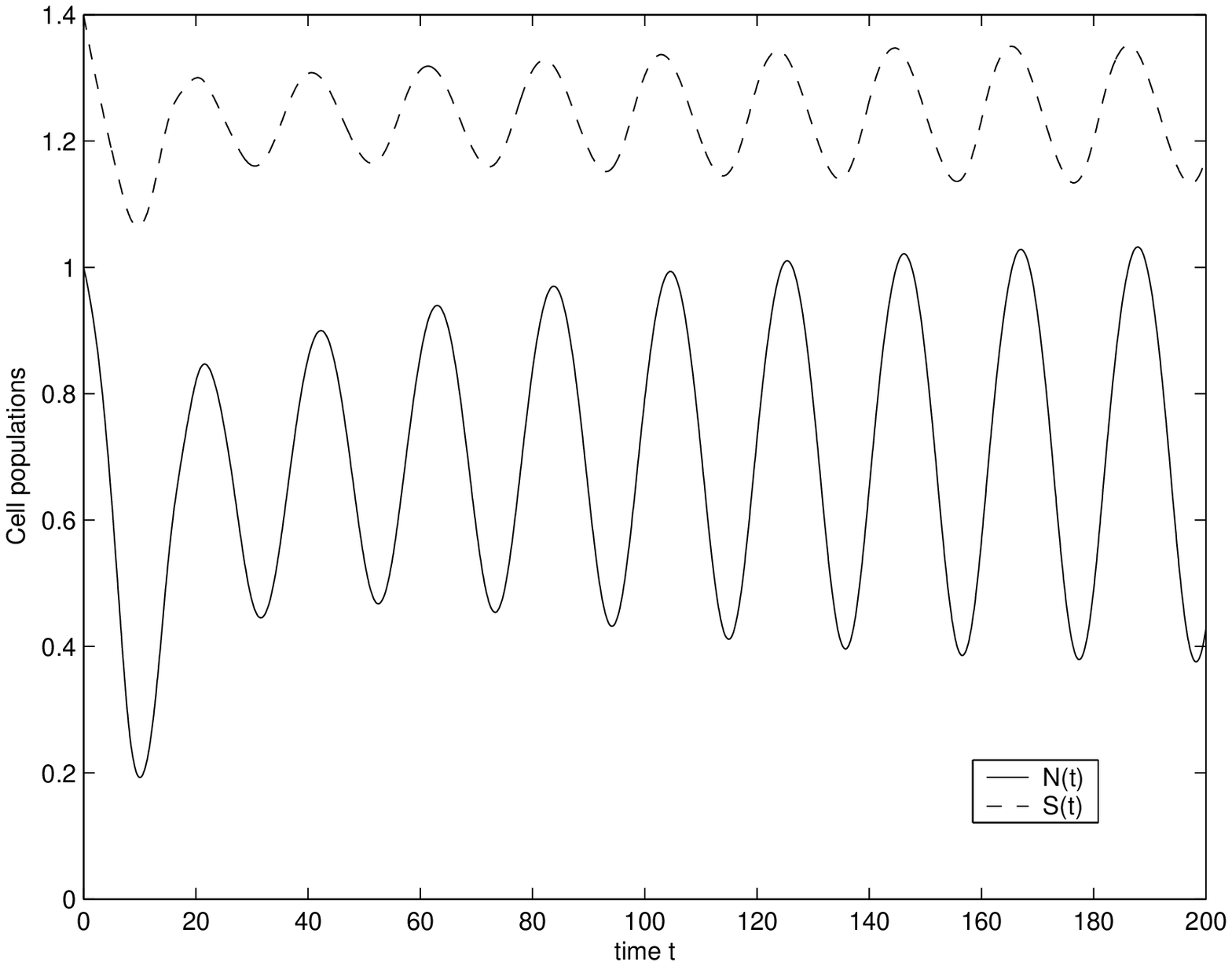}
\includegraphics[width=6cm, height=4.25cm]{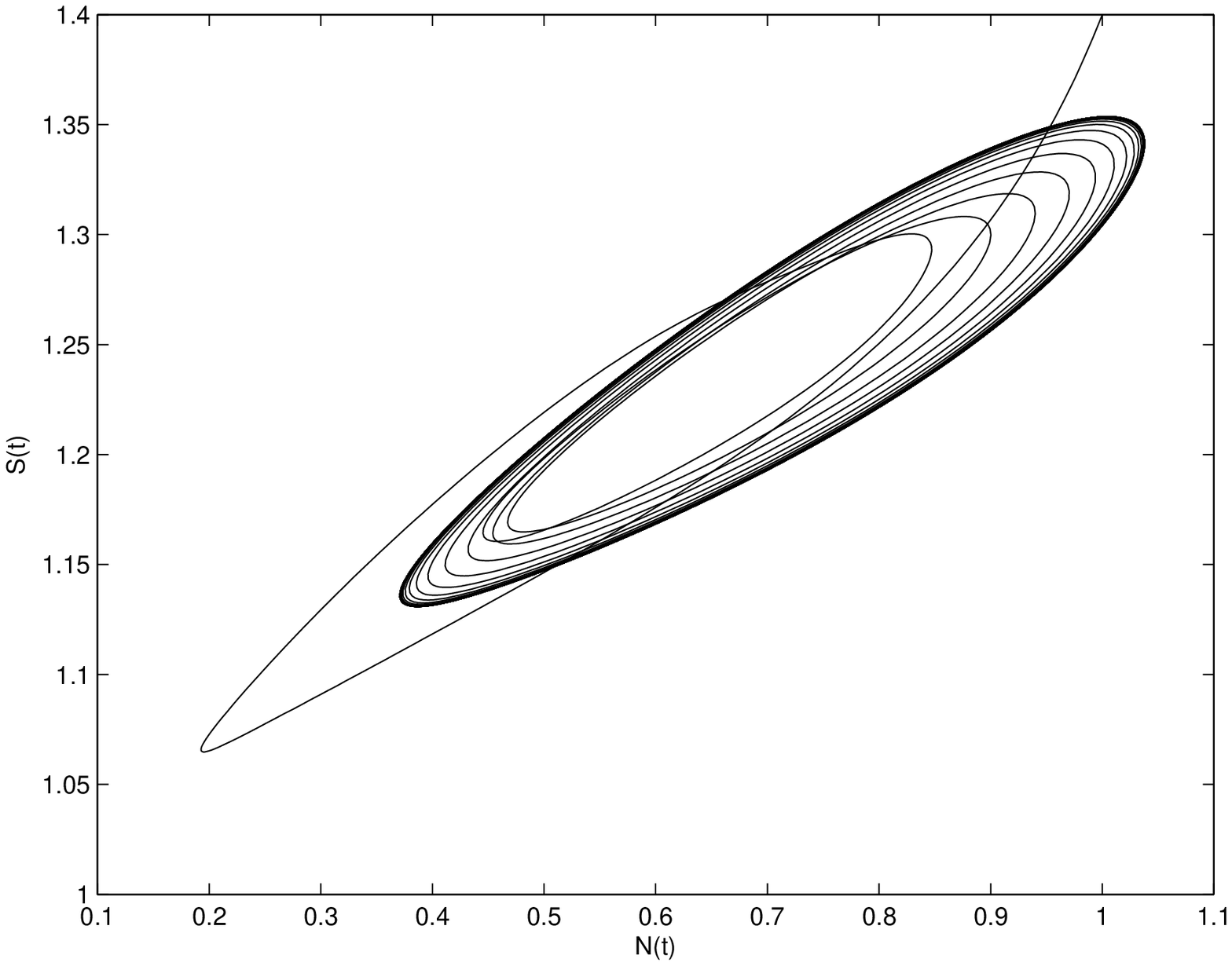}
\caption{For $\tau=7$ days, and the other parameters given by (\ref{par}) with $n=12$, long periods
oscillations are observed, with periods about 20-25 days. The steady state $E^*$ is
unstable.}\label{osc7_20}
\end{center}
\end{figure}

When $\tau$ increases after the bifurcation, one can observe oscillating solutions with longer
periods (in the order of 20 to 30 days), as it can be seen in Figure \ref{osc7_20}.

This phenomenon has already been observed by Pujo-Menjouet et al. \cite{pbm2005, pm2004}. It can be
related to diseases affecting blood cells, the so-called periodic hematological diseases
\cite{hdm1998}, which are characterized by oscillations of circulating blood cell counts with long
periods compared to the cell cycle duration. Among the wide variety of periodic hematological
diseases, we can cite chronic myelogenous leukemia \cite{acr2005, fm1999, pbm2005}, a cancer of
white blood cells with periods usually falling in the range of 70 to 80 days, and cyclical
neutropenia \cite{bbm2003, hdm1998} which is known to exhibit oscillations around 3 weeks of
circulating neutrophils (white cells), as observed on Figure \ref{osc7_20}.

\begin{figure}[pt]
\begin{center}
\includegraphics[width=6cm, height=4cm]{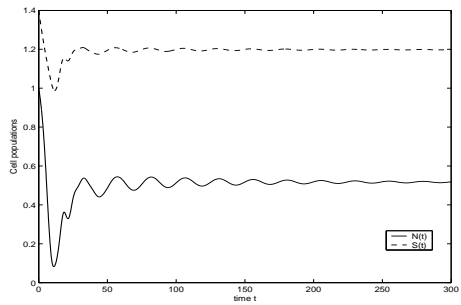}
\caption{For $\tau=9$ days, and the other parameters given by (\ref{par}) with $n=12$, damped
oscillations are observed and the steady state is stable.}\label{stab9}
\end{center}
\end{figure}

Eventually, one can note that when $\tau$ passes through the second critical value $\tau_2$,
stability switches and the steady state $(S^*,N^*)$ becomes stable again (see Figure \ref{stab9}).

\section{Discussion}

We considered a nonlinear model of blood cell dynamics in which the nonlinearity depends upon the
entire hematopoietic stem cell population, contrary to the common assumption used in previous works
\cite{acr2005_2, acr2005, bbm2003, m1978, m1979, pbm2005, pm2004} dealing with blood cell
production models. Then we were lead to the study of a new nonlinear system of two differential
equations with delay (describing the cell cycle duration) modelling the hematopoietic stem cells
dynamics.

We obtained the existence of two steady states for this model: a trivial one and a positive
delay-dependent steady state. Through sections \ref{stss} and \ref{spss}, we performed the
stability analysis of our model. We determined necessary and sufficient conditions for the global
asymptotic stability of the trivial steady state of system (\ref{eqS})--(\ref{eqN2}), which
describes the population's dying out. Using an approach proposed by Beretta and Kuang
\cite{bk2002}, we analyzed a first degree exponential polynomial characteristic equation with
delay-dependent coefficients in order to obtain the existence of a Hopf bifurcation for the
positive steady state (see Theorem \ref{theohopf}), leading to the existence of periodic solutions.

On the example presented in the previous section, we obtained long periods oscillations, which can
be related to some periodic hematological diseases (in particularly, to cyclical neutropenia
\cite{bbm2003}). This result is in keeping with previous analysis of blood cell dynamics models (as
it can be found in \cite{m1978, pbm2005, pm2004}). Periodic hematological diseases are particular
diseases mostly originated from the hematopoietic stem cell compartment. The appearance of periodic
solutions in our model with periods that can be related to the ones observed in some periodic
hematological diseases stresses the interesting properties displayed by our model. Periods of
oscillating solutions can for example be used to determine the length of cell cycles in
hematopoietic stem cell populations that cannot be directly determined experimentally.

Moreover, stability switches have been observed, due to the structure of the equations (nonlinear
equations with delay-dependent coefficients). Such a behavior had been noted in previous works
dealing with blood cell production models (see \cite{pbm2005, pm2004}), but it had never been
mathematically explained.

We can note that our assumption that proliferating and nonproliferating cells die with the same
rate may be too limitative, since Pujo-Menjouet et al. \cite{pbm2005, pm2004} already noticed that
the apoptotic rate (the proliferating phase mortality rate $\gamma$) plays an important role in the
appearance of oscillating solutions. However, by assuming that the two populations die with
different rates, we are lead to a second order exponential polynomial characteristic equation, and
the calculations are more difficult than the ones carried out in the present work. We let it for
further analysis.


\end{document}